\documentclass[letterpaper,10pt,final,openbib]{article}
\usepackage[latin1]{inputenc}
\usepackage{amsfonts}
\usepackage{mathrsfs}
\usepackage{amsmath}
\usepackage{amssymb}
\usepackage{amsthm}
\usepackage[T1]{fontenc}
\usepackage[dvips]{graphicx}

\newcommand{\ud}{\mathrm{d}}
\newcommand{\nd}{\stackrel{def}{=}}
\newcommand{\lla}{\left\langle}
\newcommand{\rra}{\right\rangle}

\newtheorem{Theorem}{Theorem}[section]
\newtheorem{Definition}[Theorem]{Definition}
\newtheorem{Proposition}[Theorem]{Proposition}
\newtheorem{Lemma}[Theorem]{Lemma}

\newtheorem{Hypothesis}[Theorem]{Hypothesis}
\newtheorem{Remark}[Theorem]{Remark}
\newtheorem{Example}[Theorem]{Example}

\begin{document}

\title{Verification theorem and construction of $\epsilon$-optimal controls
for control of abstract evolution equations}
\author{G. Fabbri\footnote {DPTEA, Universit\`a \emph{LUISS - Guido
Carli} Roma and School of Mathematics and Statistics, UNSW, Sydney
e-mail: gfabbri@luiss.it, G.Fabbri was supported by the ARC
Discovery project DP0558539.}
 \; F. Gozzi\footnote
{Dipartimento di Scienze Economiche ed Aziendali, Universit\`a
\emph{LUISS - Guido Carli} Roma, e-mail: fgozzi@luiss.it} \; and A.
\'{S}wi\c{e}ch\footnote {School of Mathematics, Georgia Institute of
Technology Atlanta, GA 30332, U.S.A., e-mail:
swiech@math.gatech.edu. A. \'{S}wi\c{e}ch was supported by NSF grant
DMS 0500270.}}
\maketitle

\begin{abstract}
\bigskip

We study several aspects of the dynamic programming approach to
optimal control of abstract evolution equations, including  a class
of semilinear partial differential equations. We introduce and prove
a verification theorem which provides a sufficient condition for
optimality. Moreover we prove sub- and superoptimality principles of
dynamic programming and give an explicit construction of
$\epsilon$-optimal controls.
\bigskip




\end{abstract}





\section{Introduction}

%

In this paper we investigate several aspects of the dynamic
programming approach to optimal control of abstract evolution
equations. The optimal control problem we have in mind has the
following form. The state equation is
\begin{equation}
\label{deterministicstateequation} \left\lbrace
\begin{array}{l}
\dot{x}(t) = Ax(t) + b(t,x(t),u(t)),\\
x(0)=x,
\end{array}
\right.
\end{equation}
$A$ is a linear, densely defined maximal dissipative operator in a
real separable Hilbert space $\mathcal{H}$, and we want to minimize
a cost functional
\begin{equation}
\label{deterministiccostfunctional} J(x;u(\cdot))= \int_0^T
L(t,x(t),u(t)) \ud t + h(x(T))
\end{equation}
over all controls
\[
u(\cdot)\in\mathcal{U}[0,T]= \{ u\colon [0,T] \to U : \; u \;
\hbox{is measurable} \},
\]
where $U$ is a metric space.

The dynamic programming approach studies the properties of the so
called value function for the problem, identifies it as a solution
of the associated Hamilton-Jacobi-Bellman (HJB) equation through the
dynamic programming principle, and then tries to use this PDE to
construct optimal feedback controls, obtain conditions for
optimality, do numerical computations, etc.. There exists an
extensive literature on the subject for optimal control of ordinary
differential equations, i.e. when the HJB equations are finite
dimensional (see for instance the books \cite{BCD, CLSW, FR, FS, Lo,
Vinter, YongZhou} and the references therein). The situation is much
more complicated for optimal control of partial differential
equations (PDE) or abstract evolution equations, i.e. when the HJB
equations are infinite dimensional, nevertheless there is by now a
large body of results on such HJB equations and the dynamic
programming approach (\cite{B1, B2, B3, B4, BaBaJe, BaDaP1, BaDaP2,
BaDaP3, BaDaP4, BaDaPPo, Ca1, CaCa, CaDaP1, CaDaP2, CaDiB, CaGoSo,
CaFr1, CaFr2, CaTe1, CaTe2, CL4, CL5, CL6, CL7, DiB, GSS, I, KoSo,
LiYong, Sh, Sri, T1, T2} and the references therein). Numerous
notions of solutions are introduced in these works, the value
functions are proved to be solutions of the dynamic programming
equations, and various verification theorems and results on
existence and explicit forms of optimal feedback controls in
particular cases are established. However, despite of these results,
so far the use of the dynamic programming approach in the resolution
of the general optimal control problems in infinite dimensions has
been rather limited. Infinite dimensionality of the state space,
unboundedness in the equations, lack of regularity of solutions, and
often complicated notions of solutions requiring the use of
sophisticated test functions are only some of the difficulties.

We will discuss two aspects of the dynamic programming approach for
a fairly general control problem: a verification theorem which gives
a sufficient condition for optimality, and the problem of
construction of $\epsilon$-optimal feedback controls.

The verification theorem we prove in this paper is an infinite
dimensional version of such a result for finite dimensional problems
obtained in \cite{Zh}. It is based on the notion of viscosity
solution (see Definitions
\ref{defdeterministicsubsol}-\ref{defdeterministicsol}). Regarding
previous result in this direction we mention \cite {CaFr1, CaFr2}
and the material in Chapter 6 \S5 of \cite{LiYong}, in particular
Theorem 5.5 there which is based on \cite {CaFr1}. We briefly
discuss this result in Remark \ref{remliyo}.

The construction of $\epsilon$-optimal controls we present here is a
fairly explicit procedure which relies on the proof of
superoptimality inequality of dynamic programming for viscosity
supersolutions of the corresponding Hamilton-Jacobi-Bellman
equation. It is a delicate generalization of such a method for the
finite dimensional case from \cite{sw}. Similar method has been used
in \cite{CLSS} to construct stabilizing feedbacks for nonlinear
systems and later in \cite{IK} for state constraint problems. The
idea here is to approximate the value function by its appropriate
inf-convolution which is more regular and satisfies a slightly
perturbed HJB inequality pointwise. One can then use this inequality
to construct $\epsilon$-optimal piecewise constant controls. This
procedure in fact gives the superoptimality inequality of dynamic
programming and the suboptimality inequality can be proved
similarly. There are other possible approaches to construction of
$\epsilon$-optimal controls. For instance under compactness
assumption on the operator $B$ (see Section 4) one can approximate
the value function by solutions of finite dimensional HJB equations
with the operator $A$ replaced by some finite dimensional operators
$A_n$ (see \cite{CL4}) and then use results of \cite{sw} directly to
construct near optimal controls. Other approximation procedures are
also possible. The method we present in this paper seems to have
some advantages: it uses only one layer of approximations, it is
very explicit and the errors in many cases can be made precise, and
it does not require any compactness of the operator $B$. It does
however require some weak continuity of the Hamiltonian and uniform
continuity of the trajectories, uniformly in $u(\cdot)$. Finally we
mention that the sub- and superoptimality inequalities of dynamic
programming are interesting on their own.

The paper is organized as follows. Definitions and the preliminary
material is presented in Section 2. Section 3 is devoted to the
verification theorem and an example where it applies in a nonsmooth
case. In Section 4 we prove sub- and superoptimality principles of
dynamic programming and show how to construct $\epsilon$-optimal
controls.

\section{Notation, definitions and background}

Throughout this paper $\mathcal{H}$ is a real separable Hilbert
space equipped with the inner product $\lla \cdot,\cdot\rra$ and the
norm $\|\cdot\|$. We recall that $A$ is a linear, densely defined
operator such that $-A$ is maximal monotone, i.e. $A$ generates a
$C_0$ semigroup of contractions $e^{sA}$, i.e.
\begin{equation}
\label{ppp1} \| e^{sA} \| \leq 1 \;\;\; \text{for all $s\geq 0$}
\end{equation}
We make the following assumptions on $b$ and $L$.
\begin{Hypothesis}
\label{hpD2onb}
\[
b\colon [0,T] \times \mathcal{H} \times U \to \mathcal{H} \;
\text{is continuous}
\]
and there exist a constant $M>0$ and a local modulus of continuity
$\omega(\cdot,\cdot)$ such that
\[
\begin{array}{ll}
\|b(t,x,u) - b(s,y,u)\| \leq M \|x-y\| + \omega(|t-s|,\|x\|\vee \|y\|)\\
\hskip 5cm \text{for all $t,s\in [0,T], \; u\in U \; x,y\in\mathcal{H}$}\\
\|b(t,0,u)\| \leq M \;\; \text{for all $(t,u)\in [0,T] \times U$}
\end{array}
\]
\end{Hypothesis}
\begin{Hypothesis}
\label{hpD3onLandh}
\[
L\colon [0,T] \times \mathcal{H} \times U \to \mathbb{R} \;\;\; and
\;\;\; h\colon \mathcal{H} \to \mathbb{R} \;\;\; \text{are
continuous}
\]
and there exist $M>0$ and a local modulus of continuity
$\omega(\cdot,\cdot)$ such that
\[
\begin{array}{ll}
|L(t,x,u) - L(s,y,u)|, \; |h(x)-h(y)| \leq \omega(\|x-y\|+|t-s|,\|x\|\vee \|y\|)\\
\hskip 5cm \text{for all $t,s\in [0,T], \; u\in U \; x,y\in\mathcal{H}$}\\
|L(t,0,u)|, |h(0)| \leq M \;\; \text{for all $(t,u)\in [0,T] \times
U$}
\end{array}
\]
\end{Hypothesis}

\begin{Remark}
Notice that if we replace $A$ and $b$ by $\tilde A=A-\omega I$ and
$b(t,x,u)$ with $\tilde b(t,x,u)= b(t,x,u) + \omega x$ the above
assumptions would cover a more general case
\begin{equation}
\label{ppp2} \| e^{sA} \| \leq e^{\omega s} \;\;\; \text{for all
$s\geq 0$}
\end{equation}
for some $\omega \geq 0$. However such $\tilde b$ does not satisfy
the assumptions of Section 4 and may not satisfy the assumptions
needed for comparison for equation (\ref{deterministicHJB}).
Alternatively, by making a change of variables $\tilde
v(t,x)=v(t,e^{\omega t}x)$ in equation (\ref{deterministicHJB}) (see
\cite{CL4}, page 275) we can always reduce the case (\ref{ppp2}) to
the case when $A$ satisfies (\ref{ppp1}).
\end{Remark}

Following the dynamic programming approach we consider a family of
problems for every $t\in[0,T], y\in \mathcal{H}$
\begin{equation}
\label{sydeterministicstate} \left\lbrace
\begin{array}{l}
\dot{x}_{t,x}(s)=A {x}_{t,x}(s) + b(s,{x}_{t,x}(s),u(s))\\
x_{t,x}(t)=x
\end{array}
\right.
\end{equation}
We will write $x(\cdot)$ for $x_{t,x}(\cdot)$ when there is no
possibility of confusion. We consider the function
\begin{equation}
\label{sydeterministiccost} J(t,x;u(\cdot))= \int_t^T L(s,x(s),u(s))
\ud t + h(x(T)),
\end{equation}
where $u(\cdot)$ is in the set of admissible controls
\[
\mathcal{U}[t,T]= \{ u\colon [t,T] \to U: \; u \hbox{ is measurable}
\}.
\]
The associated value function $V\colon [0,T]\times\mathcal{H} \to
\mathbb{R}$ is defined by
\begin{equation}
\label{deterministicvaluefunction} V(t,x)= \inf_{u(\cdot) \in
\mathcal{U}[t,T]} J(t,x;u(\cdot)).
\end{equation}
The Hamilton-Jacobi-Bellman (HJB) equation related to such optimal
control problems is
\begin{equation}
\label{deterministicHJB} \left\lbrace
\begin{array}{l}
v_t(t,x) +\lla Dv(t,x), Ax \rra + H(t,x,Dv(t,x))=0\\
v(T,x)=h(x),
\end{array}
\right.
\end{equation}
where
\[
\left\lbrace
\begin{array}{l}
H\colon [0,T]\times\mathcal{H}\times\mathcal{H} \to \mathbb{R},\\
H(t,x,p)=\inf_{u\in U} \left ( \lla p, b(t,x,u) \rra + L(t,x,u)
\right )
\end{array}
\right.
\]

The solution of the above HJB equation is understood in the
viscosity sense of Crandall and Lions \cite{CL4, CL5} which is
slightly modified here. We consider two sets of tests functions:
\[
\begin{array}{ll}
test1=\{ \varphi \in C^1((0,T)\times\mathcal{H}) \; : & \varphi
\text{ is weakly sequentially lower}\\
& \text{semicontinuous and } A^*D\varphi\in C((0,T)\times
\mathcal{H}) \}
\end{array}
\]
and
\[
\begin{array}{ll}
test2= \{ g\in C^1((0,T)\times\mathcal{H}) \; :& \exists g_0, \colon
[0,+\infty) \to [0,+\infty), \;\\
&and \; \eta\in C^1((0,T)) \text{ positive } \; s.t.\\
&g_0 \in C^1([0,+\infty)), \; g_0'(r) \geq 0 \; \forall r\geq 0, \\
& g_0'(0)=0 \; and \; g(t,x)=\eta(t)g_0(\|x\|) \\
&\forall (t,x)\in (0,T)\times \mathcal{H} \}
\end{array}
\]
We use test2 functions that are a little different from the ones
used in \cite{CL4}. The extra term $\eta(\cdot)$ in test2 functions
is added to deal with unbounded solutions. We recall that $D\varphi$
and $Dg$ stand for the Frechet derivatives of these functions.

\begin{Definition}
\label{defdeterministicsubsol} A function $v\in
C((0,T]\times\mathcal{H})$ is a (viscosity) \emph{subsolution} of
the HJB equation (\ref{deterministicHJB}) if
\[
v(T,x) \leq h(x) \;\;\; for \; all\; x\in\mathcal{H}
\]
and whenever $v-\varphi-g$ has a local maximum at $(\bar t, \bar
x)\in[0,T)\times\mathcal{H}$ for $\varphi \in test1$ and $g\in
test2$, we have
\begin{equation}
\label{eqsubsol} \varphi_t(\bar t, \bar x) + g_t(\bar t, \bar
x)+\lla A^* D \varphi(\bar t, \bar x) , \bar x \rra +H(\bar t, \bar
x, D\varphi(\bar t, \bar x)+ Dg(\bar t, \bar x)) \geq 0.
\end{equation}
\end{Definition}
\begin{Definition}
\label{defdeterministicsupersol} A function $v\in
C((0,T]\times\mathcal{H})$ is a (viscosity) \emph{supersolution} of
the HJB equation (\ref{deterministicHJB}) if
\[
v(T,x) \geq h(x) \;\;\; for \; all\; x\in\mathcal{H}
\]
and whenever $v+\varphi+g$ has a local minimum at $(\bar t, \bar
x)\in[0,T)\times\mathcal{H}$ for $\varphi \in test1$ and $g\in
test2$, we have
\begin{equation}
\label{eqsupersol} -\varphi_t(\bar t, \bar x) - g_t(\bar t, \bar x)
- \lla A^* D \varphi(\bar t, \bar x) , \bar x \rra + H(\bar t, \bar
x, -D \varphi(\bar t, \bar x)- D g(\bar t, \bar x)) \leq 0.
\end{equation}
\end{Definition}
\begin{Definition}
\label{defdeterministicsol} A function $v\in
C((0,T]\times\mathcal{H})$ is a (viscosity) \emph{solution} of the
HJB equation (\ref{deterministicHJB}) if it is at the same time a
subsolution and a supersolution.
\end{Definition}

We will be also using viscosity sub- and supersolutions in
situations where no terminal values are given in
(\ref{deterministicHJB}). We will then call a viscosity subsolution
(respectively, supersolution) simply a function that satisfies
(\ref{eqsubsol}) (respectively, (\ref{eqsupersol})).


\begin{Lemma}
\label{lemmaphi} Let Hypotheses \ref{hpD2onb} and \ref{hpD3onLandh}
hold. Let $\phi\in test1$ and $(t,x)\in(0,T)\times\mathcal{H}$. Then
the following convergence holds uniformly in $u(\cdot) \in
\mathcal{U}[t,T]$:
\begin{multline}
\lim_{s\downarrow t} \left ( \frac{1}{s-t} \left (
\varphi(s,x_{t,x}(s)) - \varphi(t,x) \right ) - \varphi_t(t,x) -
\lla A^*D\varphi (t,x),x\rra
\right.\\
\left. - \frac{1}{s-t} \int_t^s \lla D \varphi (t,x), b(t,x,u(r))
\rra \ud r \right ) =0
\end{multline}
Moreover we have for $s-t$ sufficiently small
\begin{multline}
\label{eq:explicitphi} \varphi(s,x_{t,x}(s))-\varphi(t,x) = \int_t^s
\varphi_t(r,x_{t,x}(r)) +
\lla A^*D\varphi (r,x_{t,x}(r)),x_{t,x}(r)\rra \\
+ \lla D \varphi(r,x_{t,x}(r)), b(r,x_{t,x}(r),u(r)) \rra \ud r
\end{multline}
\end{Lemma}
\begin{proof}
See \cite{LiYong} Lemma 3.3 page 240 and Proposition 5.5 page 67.
\end{proof}

\begin{Lemma}
\label{lemmag} Let Hypotheses \ref{hpD2onb} and \ref{hpD3onLandh}
hold. Let $g\in test2$ and $(t,x)\in(0,T)\times\mathcal{H}$. Then
for $s-t \to 0^+$
\begin{multline}
\label{gconv1} \frac{1}{s-t} \left ( g(s,x_{t,x}(s)) - g(t,x) \right
) \leq g_t(t,x)
\\
+ \frac{1}{s-t} \int_t^s \lla D g (t,x), b(t,x,u(r)) \rra \ud r +
o(1)
\end{multline}
where $o(1)$ is uniform in $u(\cdot) \in \mathcal{U}[t,T]$
\end{Lemma}
\begin{proof}
To prove the statement when $x \ne 0$ we use the fact that, in this
case (see \cite{LiYong} page 241, equation (3.11)),
\[
\|x_{t,x}(s)\|  \leq \|x\| + \int_t^s \lla \frac{x}{\|x\|},
b(t,x,u(r)) \rra \ud r + o(s-t)
\]
So we have
\begin{multline}
\label{eq:prooflemmag}
g(s,x_{t,x}(s))-g(t,x) = \eta(s) g_0(\|x_{t,x}(s)\|) -\eta(t)g_0(\|x\|)\\
\leq \eta(s) g_0\left ( \|x\| + \int_t^s \lla \frac{x}{\|x\|}, b(t,x,u(r)) \rra \ud r + o(s-t) \right ) - \eta(t)g_0(\|x\|) \\
\leq \eta'(t) g_0(\|x\|) (s-t) + \eta(t) g_0'(\|x\|) \left (
\int_t^s \lla \frac{x}{\|x\|}, b(t,x,u(r)) \rra
\ud r \right ) + o(s-t) \\
= g_t(t,x) (s-t) + \int_t^s \lla D g(t,x) , b(t,x,u(r)) \rra \ud r +
o(s-t)
\end{multline}
where $o(s-t)$ is uniform in $u(\cdot)$. When $x=0$, using the fact
that $g'_0(0)=0$, we get
\[
g(s,x_{t,x}(s))-g(t,x)=g_t(t,x) (s-t) + o(s-t+\|x_{t,x}(s)\|)
\]
and (\ref{gconv1}) follows upon noticing that $\|x_{t,x}(s)\|\le
C(s-t)$ for some $C$ independent of $u(\cdot) \in \mathcal{U}[t,T]$.
\end{proof}
\begin{Theorem}
\label{thexistence} Let Hypotheses \ref{hpD2onb} and
\ref{hpD3onLandh} hold. Then the value function $V$ (defined in
(\ref{deterministicvaluefunction})) is a viscosity solution of the
HJB equation (\ref{deterministicHJB}).
\end{Theorem}
\begin{proof}
The proof is quite standard and can be obtained with small changes
(due to the small differences in the definition of test2 functions)
from Theorem 2.2, page 229 of \cite{LiYong} and the proof of Theorem
3.2, page 240 of \cite{LiYong} (or from \cite{CL5}).
\end{proof}

We will need a comparison result in the proof of the verification
theorem. There are various versions of such results for equation
(\ref{deterministicHJB}) available in the literature, several
sufficient sets of hypotheses can be found in \cite{CL4, CL5}. Since
we are not interested in the comparison result itself we choose to
assume a form of comparison theorem as a hypothesis.
\begin{Hypothesis}
\label{D4deterministiccomparison} There exists a set
$\mathcal{G}\subseteq C([0,T]\times\mathcal{H})$ such that:
\begin{itemize}
\item[(i)] the value function $V$ is in $\mathcal{G}$;
\item[(ii)] if $v_1, v_2 \in \mathcal{G}$, $v_1$ is a subsolution of
the HJB equation (\ref{deterministicHJB}) and $v_2$ is a
supersolution of the HJB equation (\ref{deterministicHJB}) then
$v_1\leq v_2$.
\end{itemize}
\end{Hypothesis}
Note that from $(i)$ and $(ii)$ we know that $V$ is the only
solution of the HJB equation (\ref{deterministicHJB}) in
$\mathcal{G}$.

We will use the following lemma whose proof can be found in
\cite{YongZhou}, page 270.
\begin{Lemma}
\label{lemmaYZ} Let $g\in C([0,T];\mathbb{R})$. We extend $g$ to a
function (still denoted by $g$) on $(-\infty,+\infty)$ by setting
$g(t)=g(T)$ for $t>T$ and $g(t)=g(0)$ for $t<0$. Suppose there is a
function $\rho \in L^1(0,T;\mathbb{R})$ such that
\[
\limsup_{h\to 0^+} \frac{g(t+h) - g(t)}{h}\leq \rho(t) \;\;\; a.e.
\; t\in[0,T].
\]
Then
\[
g(\beta)-g(\alpha) \leq \int_\alpha^\beta \limsup_{h\to 0^+}
\frac{g(t+h) - g(t)}{h} \ud t\;\;\;\; \forall \;
0\leq\alpha\leq\beta\leq T.
\]
\end{Lemma}

We will denote by $B_R$ the open ball of radius $R$ centered at $0$
in $\mathcal{H}$.

\section{The verification theorem}
We first introduce a set related to a subset of the
superdifferential of a function in $C((0,T)\times\mathcal{H})$. Its
definition is suggested by the definition of a sub/super solution.
We recall that the superdifferential $D^{1,+}v(t,x)$ of $v \in
C((0,T)\times\mathcal{H})$ at $(t,x)$ is given by the pairs
$(q,p)\in \mathbb{R}\times \mathcal{H}$ such that $v(s,y) - v(t,x)
 - \left\langle p, y-x \right\rangle - q(s-t) \leq o(\|x-y\| + |t-s|)$,
and the subdifferential $D^{1,-}v(t,x)$ at $(t,x)$ is the set of all
$(q,p)\in \mathbb{R}\times \mathcal{H}$ such that $v(s,y) - v(t,x) -
\left\langle p, y-x \right\rangle - q(s-t) \geq o(\|x-y\| + |t-s|)$.
\begin{Definition}
\label{defE} Given $v\in C((0,T)\times\mathcal{H})$ and
$(t,x)\in(0,T)\times \mathcal{H}$ we define $E^{1,+} v(t,x)$ as
\[
\begin{array}{ll}
E^{1,+}v(t,x)= \{ (q,p_1,p_2)\in \mathbb{R}\times D(A^*)
\times\mathcal{H} : &
\exists \varphi\in test1, \; g\in test2\; s.t.\\
& v-\varphi-g \text{ attains a local}\\
& \text{maximum at } (t,x),\\
& \partial_t(\varphi+g)(t,x)=q,\\
& D\varphi(t,x)=p_1, \;\; Dg(t,x)=p_2\\
& and \; v(t,x)= \varphi(t,x)+g(t,x) \}
\end{array}
\]
\end{Definition}
\begin{Remark}
If we define
\[
E^{1,+}_1v(t,x) =\{(q,p)\in \mathbb{R}\times \mathcal{H} \; : \;
p=p_1+p_2 \; with \; (q,p_1,p_2) \in E^{1,+}v(t,x) \}
\]
then $E^{1,+}_1v(t,x) \subseteq D^{1,+}v(t,x)$ and in the finite
dimensional case we have \linebreak $E^{1,+}_1v(t,x) =
D^{1,+}v(t,x)$. Here we have to use $E^{1,+}v(t,x)$ instead of
$E^{1,+}_1v(t,x)$ because of the different roles of $g$ and
$\varphi$. It is not clear if the sets $E^{1,+}v(t,x)$ and
$E^{1,+}_1v(t,x)$ are convex. However if we took finite sums of
functions $\eta(t)g_0(\|x\|)$ as $test2$ functions then they would
be convex. All the results obtained are unchanged if we use the
definition of viscosity solution with this enlarged class of $test2$
functions.
\end{Remark}

\begin{Definition}
A trajectory-strategy pair $\left(x(\cdot), u(\cdot) \right)$ will
be called an {\rm admissible couple} for $(t,x)$ if $u\in{\cal
U}[t,T]$ and $x(\cdot)$ is the corresponding solution of the state
equation (\ref{sydeterministicstate}).

A trajectory-strategy pair $\left(x^*(\cdot),u^*(\cdot) \right)$
will be called an {\rm optimal couple} for $(t,x)$ if it is
admissible for $(t,x)$ and if we have
\[
-\infty < J(t,x;u^*(\cdot))\leq J(t,x;u(\cdot))
\]
for every admissible control $u(\cdot) \in {\cal U}[t,T]$.
\end{Definition}

We can now state and prove the verification theorem.

\begin{Theorem}
\label{thdeterministicverification} Let Hypotheses \ref{hpD2onb},
\ref{hpD3onLandh} and \ref{D4deterministiccomparison} hold. Let
$v\in\mathcal{G}$ be a subsolution of the HJB equation
(\ref{deterministicHJB}) such that
\begin{equation}
\label{terminalconditionvertheorem} v(T,x)=h(x) \;\;\; for\; all\;
x\; in \; \mathcal{H}.
\end{equation}

(a) We have $v(t,x) \leq V(t,x) \leq J(t,x,u(\cdot))\;\;
\forall(t,x) \in (0,T]\times\mathcal{H},
\;u(\cdot)\in\mathcal{U}[t,T]$.

(b) Let $(t,x)\in (0,T)\times H$ and let $(x_{t,x}(\cdot),
u(\cdot))$ be an admissible couple at $(t,x)$. Assume that there
exist $q\in L^1(t,T;\mathbb{R})$, $p_1\in L^1(t,T;D(A^*))$ and
$p_2\in L^1(t,T;\mathcal{H})$ such that
\begin{equation}
\label{condE} (q(s),p_1(s),p_2(s)) \in E^{1,+}v(s,x_{t,x}(s)) \;
\;\; \text{for almost all } s\in (t,T)
\end{equation}
and that
\begin{multline}
\label{condmin} \int_t^T (\lla p_1(s) + p_2(s), b(s,x_{t,x}(s),u(s))
\rra + q(s) +
\lla A^* p_1(s),x_{t,x}(s)\rra)\ud t \\
\leq \int_t^T - L(s,x_{t,x}(s),u(s)) \ud s.
\end{multline}
Then $(x_{t,x}(\cdot), u(\cdot))$ is an optimal couple at $(t,x)$
and $v(t,x)=V(t,x)$. Moreover we have equality in (\ref{condmin}).
\end{Theorem}

\begin{Remark}
It is tempting to try to prove, along the lines of Theorem 3.9,
p.243 of \cite{YongZhou}, that a condition like (\ref{condmin}) can
also be necessary if $v$ is a viscosity solution (or maybe simply a
supersolution). However this is not an easy task: the main problem
is that $E^{1,+}$ and the analogous object $E^{1,-}$ are
fundamentally different so a natural generalization of a result like
Theorem 3.9, p.243 of \cite{YongZhou} does not seem possible.
Moreover our verification theorem has some drawbacks. Condition
(\ref{condmin}) implicitly implies that $<p_2(r),Ax_{t,x}(r)>=0$
a.e. if the trajectory is in the domain of $A$. This follows from
the fact that we would then have an additional term
$<p_2(r),Ax_{t,x}(r)>$ in the integrand of the middle line of
(\ref{acbd}) so (\ref{condmin}) would also have to be an equality
with this additional term. Therefore the applicability of the
theorem is somehow limited as in practice (\ref{condmin}) may be
satisfied only if the function is ``nice" (i.e. its
superdifferential should really only consist of $p_1$). Still it
applies in some cases where other results fail (see Remarks
\ref{remliyo} and \ref{rm:controesempio}). Many issues are not fully
resolved yet and we plan to work on them in the future.
\end{Remark}

\begin{proof}
The first statement ($v\le V$) follows from Hypothesis
\ref{D4deterministiccomparison}, it remains to prove second one. The
function
\[
\left\lbrace
\begin{array}{l}
[t,T]\to \mathcal{H}\times\mathbb{R}\\
s \mapsto (b(s,x_{t,x}(s),u(s)), L(s,x_{t,x}(s),u(s))
\end{array}
\right. \]
in view of Hypotheses \ref{hpD2onb} and \ref{hpD3onLandh} is in
$L^1(t,T;\mathcal{H}\times\mathbb{R})$ (in fact it is bounded). So
the set of the right-Lebesgue points of this function that in
addition satisfy (\ref{condE}) is of full measure. We choose $r$ to
be a point in this set. We will denote $y= x_{t,x}(r)$.

Consider now two functions $\varphi^{r,y}\in test1$ and $g^{r,y}\in
test2$ such that (we will avoid the index $^{r,y}$ in the sequel)
$v\leq \varphi +g$ in a neighborhood of $(r,y)$, $v(r,y) -
\varphi(r,y) - g(r,y) =0$,$(\partial_t)(\varphi+g)(r,y))=q(r)$,
$D\phi(r,y)=p_1(r)$ and $D g(r,y)=p_2(r)$. Then for $\tau\in(r,T]$
such that $(\tau-r)$ is small enough we have by Lemmas
\ref{lemmaphi} and \ref{lemmag}
\[
\frac{v(\tau,x_{t,x}(\tau)) - v(r,y)}{\tau-r} \leq
\frac{g(\tau,x_{t,x}(\tau)) - g(r,y) }{\tau-r} +
\frac{\varphi(\tau,x_{t,x}(\tau)) - \varphi(r,y)}{\tau-r}
\]
\begin{multline}
\leq g_t(r,y) + \frac{\int_r^\tau \lla D g(r,y) , b(r,y,u(s)) \rra
\ud s}{\tau-r}
\\
+ \varphi_t(r,y) + \frac{\int_r^\tau \lla D\varphi(r,y) ,
b(r,y,u(s)) \rra \ud s}{\tau-r}+ \lla A^*D\varphi(r,y),y\rra + o(1).
\end{multline}
In view of the choice of $r$ we know that
\[
\frac{\int_r^\tau \lla D g(r,y) , b(r,y,u(s)) \rra \ud s}{\tau-r}
\xrightarrow{\tau\to r} \lla D g(r,y) , b(r,y,u(r)) \rra
\]
and
\[
\frac{\int_r^\tau \lla D \varphi(r,y) , b(r,y,u(s)) \rra \ud
s}{\tau-r} \xrightarrow{\tau\to r} \lla D \varphi(r,y) , b(r,y,u(r))
\rra.
\]
Therefore for almost every $r$ in $[t,T]$ we have
\begin{multline}
\limsup_{\tau\downarrow r} \frac{v(\tau,x_{t,x}(\tau))
- v(r,x_{t,x}(r)))}{\tau-r}\\
\leq \lla D g(r,x_{t,x}(r)) + D \varphi(r,x_{t,x}(r)),
b(r,x_{t,x}(r),u(r))\rra\\
+ g_t(r,x_{t,x}(r))+ \varphi_t(r,x_{t,x}(r)) +\lla A^* D
\varphi(r,x_{t,x}(r)),x_{t,x}(r)\rra\\
= \lla p_1(r)+p_2(r), b(r,x_{t,x}(r),u(r))\rra + q(r) + \lla A^*
p_1(r),x_{t,x}(r)\rra.
\end{multline}
We can then use Lemma \ref{lemmaYZ} and (\ref{condmin}) to obtain
\begin{multline}
\label{acbd} v(T,x_{t,x}(T)) - v(t,x)
\\
\leq \int_t^T (\lla p(r), b(r,x_{t,x}(r),u(r))\rra + q(r) + \lla A^*
 p_1(r),x_{t,x}(r)\rra)\ud r
\\
\leq \int_t^T -L(r,x_{t,x}(r),u(r)) \ud r.
\end{multline}
Thus, using (a), we finally arrive at
\begin{multline}
V(T,x_{t,x}(T)) - V(t,x) = h(x_{t,x}(T)) - V(t,x) \leq h(x_{t,x}(T))
- v(t,x)
\\
= v(T,x_{t,x}(T))-v(t,x) \leq \int_t^T -L(r,x_{t,x}(r),u(r)) \ud r
\end{multline}
which implies that $(x_{t,x}(\cdot), u(\cdot))$ is an optimal pair
and that $v(t,x)=V(t,x)$.
\end{proof}

\begin{Remark}
\label{remliyo} In the book \cite{LiYong} (page 263, Theorem 5.5)
the authors present a verification theorem (based on a previous
result of \cite{CaFr2}, see also \cite{CaFr1} for similar results)
in which it is required that the trajectory of the system remains in
the domain of $A$ a.e. for the admissible control $u(\cdot)$ in
question. This is not required here and in fact this is not
satisfied in the example of the next section.

It is shown in \cite{LiYong} (under assumptions similar to
Hypotheses \ref{hpD2onb} and \ref{hpD3onLandh}) that the couple
$x(\cdot), u(\cdot))$ is optimal if and only if
\begin{multline}
u(s) \in \bigg\lbrace u\in U \, : \, \lim_{\delta \to 0}
\frac{V((s+\delta), x(s)+ \delta(Ax(s) + b(s,x(s),u)) ) - V(s,x(s))
}{\delta} \\
= -L(s,x(s),u) \bigg\rbrace
\end{multline}
for almost every $s\in[t,T]$, where $V$ is the value function.
\end{Remark}

\subsection{An example}
We present an example of a control problem for which the value
function is a nonsmooth viscosity solution of the corresponding HJB
equation, however we can apply our verification theorem. The problem
can model a number of phenomena, for example in age-structured
population models (see \cite{Iannelli95, Anita01, Iannelli06}), in
population economics \cite{FeichtingerPrskwetzVeliov04}, optimal
technology adoption in a vintage capital context
\cite{BarucciGozzi98, BarucciGozzi01}.


Consider the state equation
\begin{equation}
\label{example1stateequation} \left\lbrace
\begin{array}{l}
\dot{x}(s) = Ax(s) + Ru(s)\\
x(t)=x
\end{array}
\right.
\end{equation}
where$A$ is a linear, densely defined maximal dissipative operator
in $\mathcal{H}$, $R$ is a continuous linear operator $R\colon
\mathbb{R}\to \mathcal{H}$, so it is of the form $R\colon u\mapsto u
\beta$ for some $\beta\in\mathcal{H}$. Let $B$ be an operator as in
Section \ref{subsuper} satisfying (\ref{bcond}). We will be using
the notation of Section \ref{subsuper}.

We will assume that $A^*$ has an eigenvalue $\lambda$ with an
eigenvector $\alpha$ belonging to the range of $B$.

We consider the functional to be minimized
\begin{equation}
\label{example1costfunctional} J(x,u(\cdot))= \int_t^T
-\left|\left\langle \alpha,x(s)\right\rangle \right| +
\frac{1}{2}u(s)^2 \ud s.
\end{equation}
We define
\[
\bar\alpha(t)\nd \int_t^T e^{(s-t)A^*} \alpha \ud s
\]
and we take $M\nd \sup_{t\in[0,T]} |\lla \bar\alpha(t), \beta
\rra|$. We consider as control set $U$ the compact subset of
$\mathbb{R}$ given by $U=[-M-1, M+1]$. So we specify the general
problem characterized by (\ref{deterministicstateequation}) and
(\ref{deterministiccostfunctional}) taking $b(t,x,u)=Ru$, $L(t,x,u)=
-\left|\left\langle \alpha,x(s)\right\rangle \right| + 1/2 u(t)^2$,
$h=0$, $U=[-M-1, M+1]$.

The HJB equation (\ref{deterministicHJB}) becomes
\begin{equation}
\label{example1HJB} \left\lbrace
\begin{array}{l}
v_t + \lla Dv, Ax \rra  -\left|\left\langle \alpha,x\right\rangle
\right|+\inf_{u\in U} \left ( \lla
u,R^*Dv\rra_{\mathbb{R}} + \frac{1}{2} u^{2} \right)=0\\
v(T,x)=0
\end{array}
\right.
\end{equation}
Note that the operator $R^*\colon \mathcal{H} \to \mathbb{R}$ can be
explicitly expressed using $\beta$ which was used to define the
operator $R$: $R^*x=\left\langle\beta,x\right\rangle$.

Now we observe that for $\left\langle \alpha,x\right\rangle<0$
(respectively $>0$) the HJB equation is the same as the one for the
optimal control problem with the objective functional $\int_t^T
\left\langle \alpha,x(s)\right\rangle + \frac{1}{2}u(s)^2 \ud s$
(respectively $\int_t^T -\left\langle \alpha,x(s)\right\rangle +
\frac{1}{2}u(s)^2 \ud s$) and it is known in the literature (see
\cite{FaggianGozzi} Theorem 5.5) that its solution is
\[
v_1(t,x)= \left\langle \bar\alpha(t),x \right\rangle - \int_t^T
\frac{1}{2} \left( R^*\bar\alpha(s)\right)^2 \ud s
\]
(respectively \[ v_2(t,x)= -\left\langle \bar\alpha(t),x
\right\rangle - \int_t^T \frac{1}{2} \left(
R^*\bar\alpha(s)\right)^2 \ud s).
\]
Note that on the separating hyperplane $\left\langle
\alpha,x\right\rangle=0$ the two functions assume the same values.
Indeed, since $\alpha$ an eigenvector for $A^*$,
\[
\bar\alpha(t) = G(t) \alpha
\]
where
$$
G(t)= \int_t^T e^{\lambda(s-t)}\ud s
$$
So, if $\left\langle \alpha,x\right\rangle=0$,
\[
\lla \bar \alpha (t) ,x \rra=0\;\;\;\; \text{ for all $t\in [0,T]$}.
\]
Therefore we can glue $v_1$ and $v_2$ writing
$$
W(t,x)=\left\{ \begin{array}{ll} v_1(t,x) & \hbox{if } \left\langle
\alpha,x\right\rangle \le 0 \\
 v_2(t,x) & \hbox{if } \left\langle
\alpha,x\right\rangle>0
 \end{array}\right.
$$
It is easy to see that $W$ is continuous and concave in $x$. We
claim that $W$ is a viscosity solution of (\ref{example1HJB}). For
$\left\langle \alpha,x\right\rangle<0$ and $\left\langle
\alpha,x\right\rangle>0$ it follows from the fact that $v_1$ and
$v_2$ are explicit regular solutions of the corresponding HJB
equations.

For the points $x$ where $\left\langle \alpha,x\right\rangle=0$ it
is not difficult to see that
\[
\left \{
\begin{array}{l}
D^{1,+} W(t,x)= \left \{\left(\frac{1}{2}\left(
R^*\bar\alpha(t)\right)^2,
\gamma G(t) \alpha \right) \; : \; \gamma \in [-1,1] \right \} \subseteq D(A^*)\\
D^{1,-} W(t,x)= \emptyset
\end{array}
\right .
\]

So we have to verify that $W$ is a subsolution on $\left\langle
\alpha,x\right\rangle=0$.
If $W - \varphi - g$ attains a maximum at $(t,x)$ with $\left\langle
\alpha,x\right\rangle=0$ we have that $p\nd(p_1 + p_2)\nd D(\varphi
+ g)(t,x)\in \left \{ \gamma G(t) \alpha \; : \; \gamma \in [-1,1]
\right \} \subseteq D(A^*)$. From the definition of test1 function
$p_1=D\varphi(t,x)\in D(A^*)$ so
$\eta(t)g_0'(|x|)\frac{x}{|x|}=p_2=Dg(t,x)\in D(A^*)$. $W(\cdot,x)$
is a $C^1$ function and then, recalling that $\lla \bar \alpha (t)
,x \rra_t=\lla G'(t) \alpha ,x \rra =0$, we have
\begin{equation}
\label{eq:exampleestimate0}
\partial_t (\varphi +g)(t,x)=\partial_t W(t,x)=\frac{1}{2}\left( R^*\bar\alpha(t)\right)^2,
\end{equation}
and for $p=\gamma \bar\alpha(t)$ we have
\begin{equation}
\label{eq:exampleestimate1} \inf_{u\in U} \left ( \lla Ru,p\rra +
\frac{1}{2} u^{2} \right) = - \frac{1}{2} \gamma^2 \left(
R^*\bar\alpha(t)\right)^2
\end{equation}
Moreover, recalling that $g_0'(|x|)\geq 0$ and $-A^*$ is monotone,
we have
\begin{multline}
\label{eq:exampleestimate2} \lla A^* p_1, x \rra = \lla A^*(p-p_2),
x \rra = \lla A^* \gamma G(t) \alpha, x \rra - \frac{g_0'(|x|)}{|x|}
\lla A^* x, x \rra \geq\\ \geq \gamma G(t) \lla A^* \alpha, x \rra =
0
\end{multline}
So, by (\ref{eq:exampleestimate0}), (\ref{eq:exampleestimate1}) and
(\ref{eq:exampleestimate2}),
\begin{multline}
\partial_t (\varphi +g)(t,x) + \lla A^* p_1 ,x\rra - \left| \left\langle \alpha,x\right\rangle\right| + \\
+\inf_{u\in U} \left ( \lla Ru,D(\varphi +g)(t,x)\rra + \frac{1}{2}
u^{2} \right) \geq \frac{1}{2} (1-\gamma^2) \left(
R^*\bar\alpha(s)\right)^2 \geq 0
\end{multline}
and so the claim in proved.

It is easy to see that both $W$ and the value function $V$ for the
problem are continuous on $[0,T]\times \mathcal{H}$ and moreover
$\psi=W$ and $\psi=V$ satisfy
\[
|\psi(t,x)-\psi(t,y)|\le C\|x-y\|_{-1} \quad\hbox{for
all}\,\,t\in[0,T], x,y\in \mathcal{H}
\]
for some $C\ge 0$. In particular $W$ and $V$ have at most linear
growth as $\|x\|\to\infty$. By Theorem \ref{thexistence}, the value
function $V$ is a a viscosity solution of the HJB equation
(\ref{example1HJB}) in $(0,T]\times\mathcal{H}$. Moreover, since
$\alpha=By$ for some $y\in \mathcal{H}$, comparison holds for
equation (\ref{example1HJB}) which yields $W=V$ on
$[0,T]\times\mathcal{H}$. (Comparison theorem can be easily obtained
by a modification of techniques of \cite{CL5} but we cannot refer to
any result there since both $V$ and $W$ are unbounded. However the
result follows directly from Theorem 3.1 together with Remark 3.3 of
\cite{Kel}. The reader can also consult the proof of Theorem 4.4 of
\cite{KeSw}. We point out that our assumptions are different from
the assumptions of the uniqueness Theorem 4.6 of \cite{LiYong}, page
250).

Therefore we have an explicit formula for the value function $V$
given by $V(t,x)=W(t,x)$. We see that $V$ is differentiable at
points $(t,x)$ if $\left\langle \alpha,x\right\rangle \ne  0$ and
\[
DV(t,x)=\left\lbrace
\begin{array}{ll}
\bar \alpha(t) & if \; \left\langle \alpha,x\right\rangle < 0\\
-\bar \alpha(t) & if \; \left\langle \alpha,x\right\rangle > 0
\end{array}
\right.
\]
and is not differentiable whenever $\left\langle
\alpha,x\right\rangle = 0$. However we can apply Theorem
\ref{thdeterministicverification} and prove the following result.
\begin{Proposition}
\label{propqui} The feedback map given by
\[
u^{op}(t,x)
= \left\lbrace
\begin{array}{ll}
- \left\langle \beta,\bar\alpha(t) \right\rangle  & if \; \left\langle \alpha,x\right\rangle \le 0\\
\left\langle \beta,\bar\alpha(t) \right\rangle & if \; \left\langle
\alpha,x\right\rangle > 0
\end{array}
\right.
\]
is optimal. Similarly, also the feedback map
\[
\bar u^{op}(t,x)
= \left\lbrace
\begin{array}{ll}
- \left\langle \beta,\bar\alpha(t) \right\rangle  & if \; \left\langle \alpha,x\right\rangle < 0\\
\left\langle \beta,\bar\alpha(t) \right\rangle & if \; \left\langle
\alpha,x\right\rangle \ge 0
\end{array}
\right.
\]
is optimal.
\end{Proposition}
\begin{proof}
Let $(t,x)\in (0,T]\times \mathcal{H}$ be the initial datum. If
$\left\langle \alpha,x\right\rangle \le 0$, taking the control
$-\left\langle \beta,\bar\alpha(t) \right\rangle$ the associated
state trajectory is
\[
x^{op}(s)= e^{(s-t)A}x - \int_t^{s}e^{(s-r)A} R(\left\langle
\beta,\bar\alpha(r) \right\rangle) \ud r
\]
and it easy to check that it satisfies $\left\langle
\alpha,x^{op}(s) \right\rangle \le 0 $ for every $s\ge t$. Indeed,
using the form of $R$ and the fact that $\alpha $ is eigenvector of
$A^*$ we get
$$
\left\langle \alpha,x^{op}(s)  \right\rangle = e^{\lambda(s-t)}
\left\langle \alpha, x \right\rangle
 - \left\langle \alpha, \beta \right\rangle
\int_t^{s}e^{\lambda(s-r)} \left\langle \beta,\bar\alpha(r)
\right\rangle \ud r
$$
$$
=e^{\lambda(s-t)}\left\langle \alpha, x \right\rangle
 - \left\langle \alpha, \beta \right\rangle^2
\int_t^{s}e^{\lambda(s-r)} G(r)  \ud r.
$$

Similarly if $\left\langle \alpha,x\right\rangle > 0$, taking the
control $\left\langle \beta,\bar\alpha(t) \right\rangle$ the
associated state trajectory is
\[
x^{op}(s)= e^{(s-t)A}x + \int_t^{s}e^{(s-r)A} R(\left\langle
\beta,\bar\alpha(r) \right\rangle) \ud r
\]
and it easy to check that it satisfies $\left\langle
\alpha,x^{op}(s) \right\rangle > 0 $ for every $s\ge t$.

We now apply Theorem \ref{thdeterministicverification} taking
$q(s)=\partial_t V(s, x^{op}(s))$,
\[
p_1(s)=\left\lbrace
\begin{array}{ll}
\bar \alpha(s) & if \; \left\langle \alpha,x^{op}(s)\right\rangle \le 0\\
-\bar \alpha(s) & if \; \left\langle \alpha,x^{op}(s)\right\rangle
>0
\end{array}
\right.
\]
and $p_2(s)=0$. It is easy to see that $(q(s), p_1(s), p_2(s)) \in
E^{1,+}V(s,x^{op}(s))$. The argument for $\bar u^{op}$ is completely
analogous.
\end{proof}

We continue by giving a specific example of the Hilbert space
$\mathcal{H}$, the operator $A$, and the data $\alpha$ and $\beta$.
This example is related  to the vintage capital problem in
economics, see e.g. \cite{BarucciGozzi01,BarucciGozzi98}. Let
$\mathcal{H}=L^2(0,1)$. Let $\{e^{tA}; \; t \ge 0\}$ be the
semigroup that, if we identify the points $0$ and $1$ of the
interval $[0,1]$, ``rotates'' the function:
\[
e^{tA}f(s) = f(t+s - [t+s])
\]
where $[\cdot]$ is the greatest natural number $n$ such that $n\leq
t+s$. The domain of $A$ will be
\[
D(A)= \left\lbrace f\in W^{1,2}(0,1) \; : \; f(0)=f(1) \right\rbrace
\]
and for all $f$ in $D(A)$ $A(f)(s) = \frac{\ud}{\ud s} f (s)$. We
choose $\alpha$ to be the constant function equal to $1$ at every
point of the interval $[0,1]$. (We can take for instance
$B=(I-\Delta)^{-\frac{1}{2}}$.) Moreover we choose
$\beta(s)=\chi_{[0,\frac{1}{2}]}(s) - \chi_{[0,\frac{1}{2}]}(s)$
($\chi_{\Omega} $ is the characteristic function of a set $\Omega$).
Consider an initial datum $(t,x)$ such that $\left\langle \alpha,
x\right\rangle =0$. In view of Proposition \ref{propqui}
an optimal strategy $u^{op}$ is
\[
u^{op}(s)=-\left\langle \beta,\bar\alpha(s) \right\rangle =0
\]
The related optimal trajectory is
\[
x^{op}(s)= e^{(s-t)A}y.
\]
\begin{Remark}\label{rm:controesempio}
We observe that, using such strategy, $\left\langle
\alpha,x^{op}(t)\right\rangle = 0$ for all $s\geq t$. So the
trajectory remains for a whole interval in a set in which the value
function is not differentiable. Anyway, applying Theorem
\ref{thdeterministicverification}, the optimality is proved.
Moreover $x$ can be chosen out of the domain of $A$ and so the
assumptions of the verification theorem given in \cite{LiYong} (page
263, Theorem 5.5) are not verified in this case.
\end{Remark}

\section{Sub- and superoptimality principles and construction of
$\epsilon$-optimal controls} \label{subsuper}

Let $B$ be a bounded linear positive self-adjoint operator on
$\mathcal{H}$ such that $A^*B$  bounded on $\mathcal{H}$ and let
$c_0\leq 0$ be a constant such that
\begin{equation}
\label{bcond} \lla (A^* B + c_0 B)x,x \rra \leq 0 \;\;\;\;\;\; for
\; all \; x\in\mathcal{H}.
\end{equation}
Such an operator always exists \cite{Renardy95} and we refer to
\cite{CL4} for various examples. Using the operator $B$ we define
for $\gamma>0$ the space $\mathcal{H}_{-\gamma}$ to be the
completion of $\mathcal{H}$ under the norm
\[
\|x\|_{-\gamma}=\|B^{\frac{\gamma}{2}}x\|.
\]
We need to impose another set of assumptions on $b$ and $L$.

\begin{Hypothesis}
\label{hp:section4} There exist a constant $K>0$ and a local modulus
of continuity $\omega(\cdot,\cdot)$ such that:
\[
\|b(t,x,u)-b(s,y,u)\| \leq K \|x-y \|_{-1} + \omega(|t-s|, \|x\|
\vee \|y\|)
\]
and
\[
|L(t,x,u)-L(s,y,u)| \leq \omega( \|x-y \|_{-1} + |t-s|, \|x\| \vee
\|y\|)
\]
\end{Hypothesis}

Let $m\geq 2$. Modifying slightly the functions introduced in
\cite{CL5} we define for a function $w:(0,T)\times \mathcal{H}\to
\mathbb{R}$ and $\epsilon,\beta,\lambda>0$ its sup- and
inf-convolutions by
\[
w^{\lambda,\epsilon,\beta}(t,x)=\sup_{(s,y)\in(0,T)\times
\mathcal{H}} \left\{w(s,y)-\frac{\|x-y\|_{-1}^2}{2\epsilon}
-\frac{(t-s)^2}{2\beta}-\lambda e^{2mK(T-s)}\|y\|^m\right\},
\]
\[
w_{\lambda,\epsilon,\beta}(t,x)=\inf_{(s,y)\in(0,T)\times
\mathcal{H}} \left\{w(s,y)+\frac{\|x-y\|_{-1}^2}{2\epsilon}
+\frac{(t-s)^2}{2\beta}+\lambda e^{2mK(T-s)}\|y\|^m\right\}.
\]

\begin{Lemma}
\label{lem2} Let $w$ be such that
\begin{equation}
\label{aaa2} w(t,x)\leq C(1+\|x\|^k)\quad(\hbox{respectively,}\,\,\,
w(t,x)\geq -C(1+\|x\|^k))
\end{equation}
on $(0,T)\times \mathcal{H}$ for some $k\geq 0$. Let $m>k$. Then:
\begin{itemize}
\item[(i)] For every $R>0$ there exists $M_{R,\epsilon,\beta}$
such that if $v=w^{\lambda,\epsilon,\beta}$ (respectively,
$v=w_{\lambda,\epsilon,\beta}$) then
\begin{equation}
\label{aaa6} |v(t,x)-v(s,y)|\leq
M_{R,\epsilon,\beta}(|t-s|+\|x-y\|_{-2})\quad
\hbox{on}\,\,\,(0,T)\times B_R
\end{equation}
\item[(ii)]
The function
\[
w^{\lambda,\epsilon,\beta}(t,x)+\frac{\|x\|_{-1}^2}{2\epsilon}
+\frac{t^2}{2\beta}
\]
is convex (respectively,
\[
w_{\lambda,\epsilon,\beta}(t,x)-\frac{\|x\|_{-1}^2}{2\epsilon}
-\frac{t^2}{2\beta}
\]
is concave).
\item[(iii)]
If $v=w^{\lambda,\epsilon,\beta}$ (respectively,
$v=w_{\lambda,\epsilon,\beta}$) and $v$ is differentiable at
$(t,x)\in (0,T)\times B_R$ then $|v_t(t,x)|\leq
M_{R,\epsilon,\beta}$, and $Dv(t,x)=Bq$, where $\|q\|\leq
M_{R,\epsilon,\beta}$
\end{itemize}
\end{Lemma}
\begin{proof}

\textbf{(i)} Consider the case $v=w^{\lambda, \epsilon, \beta}$.
Observe first that if $\|x\|\le R$ then
\begin{multline}
\label{eq:suponacompact}
w^{\lambda, \epsilon, \beta}(t,x) =\\
= \sup_{(s,y)\in(0,T)\times \mathcal{H} , \; \|y\|\leq N}
\left\{w(s,y)-\frac{\|x-y\|_{-1}^2}{2\epsilon}
-\frac{(t-s)^2}{2\beta}-\lambda e^{2mK(T-s)}\|y\|^m\right\},
\end{multline}
where $N$ depends only on $R$ and $\lambda$.

Now suppose $w^{\lambda, \epsilon, \beta}(t,x)\geq w^{\lambda,
\epsilon, \beta}(s,y)$. We choose a small $\sigma>0$ and $(\tilde t,
\tilde x)$ such that
\[
w^{\lambda, \epsilon, \beta}(t,x) \leq \sigma + w(\tilde t, \tilde
x) - \frac{\| x-\tilde x\|^2_{-1}}{2\epsilon} - \frac{(t-\tilde
t)^2}{2\beta} - \lambda e^{2mK(T-\tilde t)} \| \tilde x \|^m.
\]
Then
\begin{multline}
|w^{\lambda, \epsilon, \beta}(t,x)-w^{\lambda, \epsilon, \beta}(s,y)| \leq \sigma - \frac{\| x-\tilde x\|^2_{-1}}{2\epsilon} - \frac{(t-\tilde t)^2}{2\beta} + \frac{\|\tilde x - y \|^2_{-1}}{2\epsilon} + \frac{(\tilde t -s)^2}{2\beta} \\
\leq \sigma -  \frac{\lla B(x-y), x+y\rra}{2\epsilon} + \frac{\lla B(x-y), \tilde x \rra}{\epsilon} + \frac{(2\tilde t -t -s)(t-s)}{2\beta} \\
\leq \frac{(2R+N)}{2\epsilon} \|B(x-y)\| + \frac{2T}{2\beta} |t-s| +
\sigma
\end{multline}
and we conclude because of the arbitrariness of $\sigma$. The case
of $w_{\lambda, \epsilon, \beta}$ is similar.

\textbf{(ii)} It is a standard fact, see for example the Appendix of
\cite{Userguide}.

\textbf{(iii)} The fact that $|v_t(t,x)|\leq M_{R,\epsilon,\beta}$
is obvious. Moreover  if $\alpha>0$ is small and $\|y\|=1$ then
\[
\alpha M_{R,\epsilon,\beta}\|y\|_{-2}\geq |v(t,x+\alpha y)-v(x)|=
\alpha |\lla Dv(t,x),y\rra|+o(\alpha)
\]
which upon dividing by $\alpha$ and letting $\alpha\to 0$ gives
\[
|\lla Dv(t,x),y\rra|\leq M_{R,\epsilon,\beta}\|y\|_{-2}
\]
which then holds for every $y\in \mathcal{H}$. This implies that
$\lla Dv(t,x),y\rra$ is a bounded linear functional in
$\mathcal{H}_{-2}$ and so $Dv(t,x)=Bq$ for some $q\in \mathcal{H}$.
Since $|\lla q,By\rra|\leq M_{R,\epsilon,\beta}\|By\|$ we obtain
$\|q\|\leq M_{R,\epsilon,\beta}$.
\end{proof}

\begin{Lemma}
\label{lem1} Let Hypotheses \ref{hpD2onb}, \ref{hpD3onLandh} and
\ref{hp:section4} be satisfied. Let $w$ be a locally bounded
viscosity subsolution (respectively, supersolution) of
(\ref{deterministicHJB}) satisfying (\ref{aaa2}). Let $m>k$. Then
for every $R,\delta>0$ there exists a non-negative function
$\gamma_{R,\delta}(\lambda,\epsilon,\beta)$, where
\begin{equation}
\label{aaa3} \lim_{\lambda\to 0}\limsup_{\epsilon\to
0}\limsup_{\beta\to 0} \gamma_{R,\delta}(\lambda,\epsilon,\beta)=0,
\end{equation}
such that $w^{\lambda,\epsilon,\beta}$ (respectively,
$w_{\lambda,\epsilon,\beta}$) is a viscosity subsolution
(respectively, supersolution) of
\begin{equation}
\label{aaa4} v_t(t,x) +\lla Dv(t,x), Ax \rra + H(t,x,Dv(t,x))=
-\gamma_{R,\delta}(\lambda,\epsilon,\beta)\quad\hbox{in}\,\,\,
(\delta,T-\delta)\times B_R
\end{equation}
(respectively,
\begin{equation}
\label{aaa5} v_t(t,x) +\lla Dv(t,x), Ax \rra + H(t,x,Dv(t,x))=
\gamma_{R,\delta}(\lambda,\epsilon,\beta)\quad\hbox{in}\,\,\,
(\delta,T-\delta)\times B_R)
\end{equation}
for $\beta$ sufficiently small (depending on $\delta$).
\end{Lemma}

\begin{proof}
The proof is similar to the proof of Proposition 5.3 of \cite{CL5}.
We notice that $w^{\lambda,\epsilon,\beta}$ is bounded from above.

Let $(t_0,x_0)\in (\delta, T-\delta)\times \mathcal{H}$ be a local
maximum of $w^{\lambda,\epsilon,\beta}-\phi-g$. We can assume that
the maximum is global and strict (see Proposition 2.4 of \cite{CL5})
and that $w^{\lambda,\epsilon,\beta}-\phi-g\to -\infty$ as
$\|x\|\to\infty$ uniformly in $t$. In view of these facts and
(\ref{eq:suponacompact}) we can choose $S>2\|x_0\|$, depending on
$\lambda$ such that, for all $\|x\|+\|y\| >S-1$ and $s,t\in (0,T)$,
\begin{multline}
\label{eq:lessthen-1} w(s,y) - \frac{1}{2\epsilon} \|(x-y)\|_{-1}^2
- \frac{(t-s)^2}{2\beta} - \lambda e^{2mK(T-s)}\|y\|^m - \phi(t,x) -
g(t,x)
\\
\leq w(t_0,x_0) - \lambda e^{2mK(T-t_0)}\|x_0\|^m -\phi(t_0,x_0) -
g(t_0,x_0) -1.
\end{multline}
We can then use a perturbed optimization technique of \cite{CL5}
(see page 424 there) which is a version of the Ekeland-Lebourg Lemma
\cite{EkelandLebourg77} to obtain for every $\alpha>0$ elements
$p,q\in\mathcal{H}$ and $a,b\in\mathbb{R}$ with $\|p\|, \|q\|\leq
\alpha$ and $|a|,|b|\leq \alpha$ such that the function
\begin{multline}
\psi(t,x,s,y) \nd w(s,y) - \frac{1}{2\epsilon} \|(x-y)\|_{-1}^2 -
\frac{(t-s)^2}{2\beta} - \lambda e^{2mK(T-s)}\|y\|^m
\\
-g(t,x)-\phi(t,x) - \lla Bp,y \rra -\lla Bq, x \rra - at - bs
\end{multline}
attains a local maximum $(\bar t, \bar x, \bar s, \bar y)$ over
$[\delta/2,T-\delta/2]\times B_S \times [\delta/2,T-\delta/2]\times
B_S$. It follows from (\ref{eq:lessthen-1}) that if $\alpha$ is
sufficiently small then $\|\bar x\|, \|\bar y\| \leq S-1$.

By possibly making $S$ bigger we can assume that $(0,T)\times B_S$
contains a maximizing sequence for
\[
\sup_{(s,y)\in(0,T), \; \|y\|\leq N}
\left\{w(s,y)-\frac{\|x_0-y\|_{-1}^2}{2\epsilon}
-\frac{(t_0-s)^2}{2\beta}-\lambda e^{2mK(T-s)}\|y\|^m\right\}.
\]
Then
\[
\psi(\bar t, \bar x, \bar s, \bar y) \geq w^{\lambda,\epsilon,\beta}
(t_0, x_0) - \phi(t_0, x_0) - g(t_0, x_0) -C\alpha
\]
where the constant $C$ does not depend on $\alpha >0$, and
\[
\psi(\bar t, \bar x, \bar s, \bar y) \leq w^{\lambda,\epsilon,\beta}
(\bar t, \bar x) - \phi(\bar t, \bar x) - g(\bar t, \bar x) +
C\alpha.
\]
Therefore, since $(t_0, x_0)$ is a strict maximum, we have that
$(\bar t, \bar x)\xrightarrow{\alpha\downarrow 0} (t_0, x_0)$ and so
for small $\alpha$ $\bar t\in(\delta , T-\delta)$. It then easily
follows that if $\beta$ is big enough (depending on $\lambda$ and
$\delta$) then $\bar s \in (\delta/2,T-\delta/2)$.

Moreover, standard arguments (see for instance \cite{I}) give us

\begin{equation}
\label{eq:stimasuepsilon} \lim_{\beta\to 0}\limsup_{\alpha\to 0}
\frac{|\bar s -\bar t|^2}{2\beta} =0,
\end{equation}
\begin{equation}
\label{eq:stimasuepsilon1} \lim_{\epsilon\to 0}\limsup_{\beta\to
0}\limsup_{\alpha\to 0} \frac{|\bar x -\bar y|^2_{-1}}{2\epsilon}=0.
\end{equation}

We can now use the fact that $w$ is a subsolution to obtain
\begin{multline}
-\frac{(\bar t-\bar s)}{\beta} - 2\lambda mKe^{2mK(T- \bar s)}\|\bar y\|^m + b - \frac{\lla A^*B(\bar x - \bar y), \bar y \rra}{\epsilon} + \lla A^*Bp, \bar y \rra \\
+ H \left (\bar s, \bar y, \frac{1}{\epsilon}B(\bar y - \bar x) +
\lambda m e^{2mK(T-\bar s)} \|y\|^{m-1} \frac{y}{\|y\|} + Bp \right
) \geq 0.
\end{multline}
We notice that
\[
-\frac{(\bar t-\bar s)}{\beta}=\phi_t(\bar t , \bar x) + g_t(\bar t,
\bar x) +a
\]
and
\[
\frac{1}{\epsilon}B(\bar y - \bar x) = D\phi (\bar t, \bar x) +
Dg(\bar t, \bar x) + Bq
\]
which in particular implies that $Dg(\bar t, \bar x)\in D(A^*)$,
i.e. $\bar x\in D(A^*)$, and so it follows that $\lla A^*\bar
x,Dg(\bar t, \bar x)\rra \le 0$. Therefore using this, the
assumptions on $b$ and $L$, and (\ref{eq:stimasuepsilon}) and
(\ref{eq:stimasuepsilon1}) we have
\begin{multline}
\phi_t(\bar t , \bar x) + g_t(\bar t, \bar x) + \lla \bar x, A^*D\phi(\bar t,\bar x)\rra + H \left (\bar t, \bar x, D\phi (\bar t, \bar x) + Dg(\bar t, \bar x) \right )  \\
\geq 2\lambda mKe^{2mK(T- \bar s)}\|\bar y\|^m - \lla A^*Bp, \bar y \rra -a -b\\
- \lla (\bar y - \bar x), A^*\frac{1}{\epsilon}B(\bar y - \bar x)\rra - \lla \bar x, A^*Dg(\bar t, \bar x) + A^*Bq)\rra \\
+ H \left (\bar t, \bar x, \frac{1}{\epsilon}B(\bar y - \bar x) -Bq \right ) - H \left (\bar s, \bar y, \frac{1}{\epsilon}B(\bar y - \bar x) + \lambda m e^{2mK(T-\bar s)} \|y\|^{m-1} \frac{y}{\|y\|} \right ) \\
\geq 2\lambda mKe^{2mK(T- \bar s)}\|\bar y\|^m
-C_{\lambda,\epsilon}\alpha
+\frac{c_0}{\epsilon} \|\bar x - \bar y \|_{-1}^2 \\
-K\|\bar x - \bar y \|_{-1}\frac{\|B(\bar x - \bar y) \|}{\epsilon}
-\gamma_{\lambda,\epsilon}(|\bar t-\bar s|)
-\lambda m(M+K\|\bar y\|)e^{2mK(T- \bar s)}\|\bar y\|^{m-1}\\
\geq -C_{\lambda,\epsilon}\alpha
-\gamma(\lambda,\epsilon,\beta,\alpha)
\end{multline}
for some $\gamma(\lambda,\epsilon,\beta,\alpha)$ such that
\[
\lim_{\lambda\to 0}\limsup_{\epsilon\to 0}\limsup_{\beta\to 0}
\limsup_{\alpha\to 0} \gamma(\lambda,\epsilon,\beta,\alpha)=0.
\]
We obtain the claim by letting $\alpha\to 0$. The proof for
$w_{\lambda,\beta,\epsilon}$ is similar.
\end{proof}

\begin{Remark}
Similar argument would also work for problems with discounting if
$w$ was uniformly continuous in $|\cdot|\times\|\cdot\|_{-1}$ norm
uniformly on bounded sets of $(0,T)\times\mathcal{H}$. Moreover in
some cases the function $\gamma_{R,\delta}$ could be explicitly
computed. For instance if $w$ is bounded and
\begin{equation}
\label{aaa12} |w(t,x)-w(s,y)|\leq
\sigma(\|x-y\|_{-1})+\sigma_1(|t-s|;\|x\|\vee\|y\|)
\end{equation}
for $t,s\in(0,T), \|x\|,\|y\|\in \mathcal{H}$, we can replace
$\lambda e^{2mK(T- \bar s)}\|\bar y\|^m$ by $\lambda\mu(y)$ for some
radial nondecreasing function $\mu$ such that $D\mu$ is bounded and
$\mu(y)\to+\infty$ as $\|y\|\to\infty$ (see \cite{CL5}, page 446).
If we then replace the order in which we pass to the limits we can
get an explicit (but complicated) form for $\gamma_{R,\delta}$
satisfying
\[
\lim_{\epsilon\to 0}\limsup_{\lambda\to 0}\limsup_{\beta\to 0}
\gamma_{R,\delta}(\epsilon,\lambda,\beta)=0.
\]
The proofs of Theorem 3.7 and Proposition 5.3 in \cite{CL5} can give
hints how to do this.

\end{Remark}

\begin{Lemma}
\label{lem3} Let the assumptions of Lemma \ref{lem1} be satisfied.
Then:
\begin{itemize}
\item[(a)]
If $(a,p)\in D^{1,-}w^{\lambda,\epsilon,\beta}(t,x)$ for $(t,x)\in
(\delta,T-\delta)\times B_R$ then
\begin{equation}
\label{eq:lem3a} a +\lla A^*p, x \rra + H(t,x,p) \geq
-\gamma_{R,\delta}(\lambda,\epsilon,\beta)
\end{equation}
for $\beta$ sufficiently small.
\item[(b)]
If in addition $H(s,y,q)$ is weakly lower-semicontinuous with
respect to the $q$-variable and $(a,p)\in
D^{1,+}w_{\lambda,\epsilon,\beta}(t,x)$ for $(t,x)\in
(\delta,T-\delta)\times B_R$ is such that
$Dw_{\lambda,\epsilon,\beta}(t_n,x_n)\rightharpoonup p$ for some
$(t_n,x_n)\to (t,x)$, $(t_n,x_n)\in(\delta,T-\delta)\times B_R$,
then
\[
a +\lla A^*p, x \rra + H(t,x,p) \leq
\gamma_{R,\delta}(\lambda,\epsilon,\beta)
\]
for $\beta$ sufficiently small.
\end{itemize}
\end{Lemma}

\begin{Remark}
The Hamiltonian $H$ is weakly lower-semicontinuous with respect to
the $q$-variable for instance if $U$ is compact. To see this we
observe that thanks to the compactness of $U$ the infimum in the
definition of the Hamiltonian is a minimum. Let now
$q_n\rightharpoonup q$ and let
\[
H(s,y,q_n)= \lla q_n, b(s, y, u_n)\rra + L(s, y, u_n)
\]
for some $u_n\in U$. Passing to a subsequence if necessary we can
assume that $u_n\longrightarrow \bar u$, and then passing to the
limit in the above expression we obtain
\[
\liminf_{n\to\infty}H(s,y,q_n)= \lla q, b(s, y, \bar u)\rra + L(s,
y,\bar u) \geq H(s,y,q).
\]
We also remark that since $H$ is concave in $q$ it is weakly
upper-semicontinuous in $q$. Therefore in (b) the Hamiltonian $H$ is
assumed to be weakly continuous in $q$.
\end{Remark}

\begin{proof} {\it (of Lemma \ref{lem3})}
Recall first that for a convex/concave function $v$ its
sub/super-differential at a point $(s,z)$ is equal to
\[
\overline{\hbox{conv}}\{((a,p):v_t(s_n,z_n)\to a,
Dv(s_n,z_n)\rightharpoonup p, s_n\to s, z_n\to z\}
\]
(see \cite{Pr}, page 319).

\textbf{(a)} \textbf{Step 1}: Denote $v=w^{\lambda,\epsilon,\beta}$.
At points of differentiability, it follows from Lemma
\ref{lem2}(iii) and the ``semiconvexity" (see Lemma \ref{lem2}(ii)) of
$w^{\lambda,\epsilon,\beta}$ that there exists a test1 function
$\varphi$ such that $v-\varphi$ has a local maximum and the result
then follows from Lemma \ref{lem1}.

\textbf{Step 2}: Consider first the case $Dv(t_n,x_n)\rightharpoonup
p$ with  $(t_n,x_n)\to (t,x)$. From Lemma \ref{lem2} (iii)
$Dv(t_n,x_n)=Bq_n$ with $\|q_n\|\leq M_{R,\epsilon,\beta}$, so, it
is always possible to extract a subsequence $q_{n_k}\rightharpoonup
q$ for some $q\in \mathcal{H}$. Then
$Dv(t_{n_k},x_{n_k})=Bq_{n_k}\rightharpoonup Bq$ and $Bq=p$.
Therefore
\[
\lla A^*B q_{n_k},x_{n_k} \rra = \lla q_{n_k},(A^*B)^*x_{n_k} \rra
\longrightarrow \lla q,(A^*B)^*x \rra = \lla A^*B q,x \rra = \lla
A^*p,x \rra
\]
Moreover, since $H$ is concave in $p$ it is weakly
upper-semicontinuous so we have
\[
H(t,x,p)\geq \limsup_{k\to +\infty} H(t_{n_k},x_{n_k},Dv(t_{n_k},x_{n_k}))
\]
and we conclude from Step 1.

\textbf{Step 3}: If $p$ is a generic point of
$\overline{\hbox{conv}}\{p:Dv(t_n,x_n)\rightharpoonup p,
(t_n,x_n)\to (t,x)\}$, i.e. $p=\lim_{n\to\infty}
\sum_{i=1}^n\lambda_i^nBq_i^n$, where $\sum_{i=1}^n\lambda_i^n=1,
\|q_i^n\| \leq M_{R,\epsilon,\beta}$, and the $Bq_i^n$ are weak
limits of gradients. By passing to a subsequence if necessary we can
assume that $\sum_{i=1}^n\lambda_i^nq_i^n\rightharpoonup q$ and
$p=Bq$. But then
\[
\lla A^*\left(\sum_{i=1}^n\lambda_i^nBq_i^n\right), x_n\rra =\lla
A^*B\left(\sum_{i=1}^n\lambda_i^nq_i^n\right), x_n\rra \to \lla
A^*Bq, x \rra=\lla A^*p, x \rra
\]
as $n\to\infty$. The result now follows from Step 2 and the
concavity of
\[
p\mapsto \lla A^*p, x \rra + H(t,x,p).
\]

\textbf{(b)} As in (a) at the points of differentiability the claim
follows from Lemmas \ref{lem2} and \ref{lem1}. Denote $v=w_{\lambda,
\epsilon, \beta}$. If $Dv(t_n, x_n)\rightharpoonup p$ for some
$(t_n,x_n)\to (t,x)$, $(t_n,x_n)\in(\delta,T-\delta)\times B_R$ we
have that
\begin{equation}
\label{aabba} v_t(t_n, x_n)+\lla A^*Dv(t_n, x_n), x_n\rra +
H(t_n,x_n,Dv(t_n, x_n)) \leq
\gamma_{R,\delta}(\lambda,\epsilon,\beta).
\end{equation}
Observing as in Step 2 of (a) that
\[
\lla A^*Dv(t_n, x_n), x_n\rra \to \lla A^* p ,x \rra
\]
we can pass to the limit in (\ref{aabba}), using the weak
lower semicontinuity of $H$ with respect to the third variable, to get
\[
a+ \lla A^* p ,x \rra + H(t, x, p) \leq
\gamma_{R,\delta}(\lambda,\epsilon,\beta).
\]
\end{proof}

\begin{Theorem}
\label{th1} Let the assumptions of Lemma \ref{lem1} be satisfied and
let $w$ be a function such that for every $R>0$ there exists a
modulus $\sigma_R$ such that
\begin{equation}
\label{aaa1} |w(t,x)-w(s,y)|\leq
\sigma_R(|t-s|+\|x-y\|_{-1})\quad\hbox{for}\,\,\, t,s\in(0,T),
\|x\|,\|y\|\leq R.
\end{equation}
Then:
\begin{itemize}

\item[(a)]
If $w$ is a viscosity subsolution of (\ref{deterministicHJB})
satisfying (\ref{aaa2}) for subsolutions then for every $0<t<t+h<T$,
$x\in \mathcal{H}$
\begin{equation}
\label{aaa7} w(t,x)\leq \inf_{u(\cdot) \in \mathcal{U}[t,T]}
\left\{\int_t^{t+h} L(s,x(s),u(s)) \ud s +w(t+h,x(t+h))\right\}.
\end{equation}
\item[(b)]
Assume in addition that $H(s,y,q)$ is weakly lower-semicontinuous in
$q$ and that for every $(t,x)$ there exists a modulus $\omega_{t,x}$
such that
\begin{equation}
\label{aaa8} \|{x}_{t,x}(s_2)-{x}_{t,x}(s_1)\|\leq
\omega_{t,x}(s_2-s_1)
\end{equation}
for all $t\leq s_1\leq s_2\leq T$ and all $u(\cdot)\in
\mathcal{U}[t,T]$, where ${x}_{t,x}(\cdot)$ is the solution of
(\ref{sydeterministicstate}). If $w$ is a viscosity supersolution of
(\ref{deterministicHJB}) satisfying (\ref{aaa2}) for supersolutions
then for every $0<t<t+h<T, x\in H$, and $\nu>0$ there exists a
piecewise constant control $u_{\nu}\in\mathcal{U}[t,T]$ such that
\begin{equation}
\label{aaa9} w(t,x)\geq \int_t^{t+h} L(s,x(s),u_{\nu}(s)) \ud s
+w(t+h,x(t+h))-\nu.
\end{equation}
In particular we obtain the superoptimality principle
\begin{equation}
\label{aaa10} w(t,x)\geq \inf_{u(\cdot) \in \mathcal{U}[t,T]}
\left\{\int_t^{t+h} L(s,x(s),u(s)) \ud s +w(t+h,x(t+h))\right\}
\end{equation}
and if $w$ is the value function $V$ we have existence (together
with the explicit construction) of piecewise constant $\nu$-optimal
controls.
\end{itemize}
\end{Theorem}
\begin{proof}
We will only prove $(b)$ as the proof of $(a)$ follows the same
strategy after we fix any control $u(\cdot)$ and is in fact much
easier. We follow the ideas of \cite{sw} (that treats the finite
dimensional case).

\textbf{Step 1}. Let $n\geq 1$. We approximate $w$ by
$w_{\lambda,\epsilon,\beta}$ with $m>k$. We notice that for any
$u(\cdot)$ if $x_{t,x}(\cdot)$ is the solution of
(\ref{sydeterministicstate}) then
\[
\sup_{t\le s\le T}\|x_{t,x}(s)\|\le R=R(T,\|x\|).
\]

\textbf{Step 2}. Take any $(a,p)\in
D^{1,+}w_{\lambda,\epsilon,\beta}(t,x)$ as in Lemma \ref{lem3}$(b)$
(i.e. $p$ is the weak limit of derivatives nearby). Such elements
always exist because $w_{\lambda,\epsilon,\beta}$ is ``semiconcave".
Then we choose $u_1\in U$ such that
\begin{equation}
\label{aaa11} a +\lla A^*p, x \rra + \lla p, b(t,x,u_1) \rra +
L(t,x,u_1) \leq
\gamma_{R,\delta}(\lambda,\epsilon,\beta)+\frac{1}{n^2}.
\end{equation}
By the ``semiconcavity" of $w_{\lambda,\epsilon,\beta}$ 
\begin{equation}
\label{aaa121} w_{\lambda,\epsilon,\beta}(s,y)\leq
w_{\lambda,\epsilon,\beta}(t,x) +a(s-t)+\lla p, y-x \rra
+\frac{\|x-y\|_{-1}^2}{2\epsilon} +\frac{(t-s)^2}{2\beta}.
\end{equation}
But the right hand side of the above inequality is a test1 function
so if $s\geq t$ and $x(s)= x_{t,x}(s)$ with constant control
$u(s)=u_1$, we can use (\ref{eq:explicitphi}) and write
\begin{multline}
\label{aaa13} \bigg | \frac{a(s-t) + \lla p, x(s) - x \rra +
\frac{\| x(s) - x \|_{-1}^2}{2\epsilon} + \frac{(s-t)^2}{2\beta}}{s-t}
\\  - \left ( a+ \lla p, b(t,x,u_1) \rra + \lla A^* p , x \rra
\right ) \bigg |
\\
\leq\frac{|t-s|}{2\beta}
+
\left | \frac{\int_t^s \lla A^*p, x(r) - x \rra \ud r}{s-t} \right | 
\\ 
+ \left| \frac{\int_t^s \lla p, b(r,x(r),u_1) - b(t,x,u_1)
\rra \ud r}{s-t} \right | + \left | \frac{\int_t^s \lla A^*B(x(r) -x), x(r) \rra \ud r}{\epsilon(s-t)} \right | 
\\
+  \left | \frac{\int_t^s \lla B(x(r) -x), b(r,x(r),u_1) \rra \ud
r}
{\epsilon(s-t)} \right | \\
\leq \omega_{t,x}' (|s-t| + \sup_{t\le r\le s}\|x(r) - x\|) \leq
\tilde\omega_{t,x}(s-t)
\end{multline}
for some moduli $\omega'_{t,x}$ and $\tilde\omega_{t,x}$ that depend
on $(t,x),\epsilon,\beta$ but not on $u_1$. We can now use (\ref{aaa11}),
(\ref{aaa121}) and (\ref{aaa13}) to estimate
\begin{multline}
\label{wzwza}
\frac{w_{\lambda,\epsilon,\beta}(t+\frac{h}{n},x(t+\frac{h}{n}))-
w_{\lambda,\epsilon,\beta}(t,x)}{h/n}
\\
\leq \tilde\omega_{t,x}\left(\frac{h}{n}\right) +
\gamma_{R,\delta}(\lambda,\epsilon,\beta)+\frac{1}{n^2}- L(t,x,u_1)
\end{multline}

\textbf{Step 3}. Denote $t_i=t+\frac{(t-1)h}{n}$ for $i=1,...,n$. We
now repeat the above procedure starting at $x(t_2)$ to abtain $u_2$
satisfying (\ref{wzwza}) with $(t_2,x(t_2))$ replaced by
$(t_3,x(t_3))$, $(t,x)=(t_1,x(t_1))$ replaced by $(t_2,x(t_2))$, and
$u_1$ replaced by $u_2$. After $n$ iterations of this process we
obtain a piecewise constant control $u^{(n)},$ where
$u^{(n)}(s)=u_i$ if $s\in [t_i,t_{i+1})$. Then if $x(r)$ solves
(\ref{sydeterministicstate}) with the control $u^{(n)}$ we have

\begin{multline}
\frac{w_{\lambda,\epsilon,\beta}(t+{h},x(t+{h}))-
w_{\lambda,\epsilon,\beta}(t,x)}{h/n}\\
\leq \tilde\omega_{t,x}\left(\frac{h}{n}\right)n +
\gamma_{R,\delta}(\lambda,\epsilon,\beta)n+\frac{n}{n^2}-
\sum_{i=1}^n L(t_{i-1},x(t_{i-1}),u_i).\nonumber
\end{multline}
We remind that (\ref{aaa8}) is needed here to guarantee that
$\sup_{t_{i-1}\le r\le t_i}\|x(r)-x(t_{i-1})\|$ is independent of
$u_i$ and $x(t_{i-1})$ and depends only on $x$ and $t$. We then
easily obtain
\begin{multline}
{w_{\lambda,\epsilon,\beta}(t+{h},x(t+{h}))-
w_{\lambda,\epsilon,\beta}(t,x)}  \\
\leq \tilde\omega_{t,x}\left(\frac{h}{n}\right)h +
\gamma_{R,\delta}(\lambda,\epsilon,\beta)h +\frac{h}{n^2}-
\int_t^{t+h} L(r,x(r),u^{(n)}) \ud r + \tilde\omega'_{t,x}
\left(\frac{h}{n}\right)h
\end{multline}
for some modulus $\tilde\omega'_{t,x}$, where we have used
Hypothesis \ref{hp:section4} and (\ref{aaa8}) to estimate how the
sum converges to the integral. We now finally notice that it follows
from (\ref{aaa1}) that
\[
|w_{\lambda,\epsilon,\beta}(s,y)-w(s,y)|\leq \tilde\sigma_R(\lambda
+ \epsilon + \beta;R) \quad\hbox{for}\,\,\,
s\in(\delta,T-\delta),\|y\|\leq R,
\]
where the modulus $\tilde\sigma_R$ can be explicitly calculated from
$\sigma_R$. Therefore, choosing $\beta,\lambda, \epsilon$ small and
then $n$ big enough, and using (\ref{aaa3}), we arrive at (\ref{aaa9}). 
\end{proof}

We show below one example when condition (\ref{aaa8}) is satisfied.
\begin{Example}
Condition (\ref{aaa8}) holds for example if $A=A^*$, it generates a
differentiable semigroup, and $\|Ae^{tA}\|\le C/t^\delta$ for some
$\delta<2$. Indeed under these assumptions, if $u(\cdot)\in
\mathcal{U}[t,T]$ and writing $x(s)=x_{t,x}(s)$, we have
\[
\|(A+I)^{\frac{1}{2}}x(s)\|\leq \|(A+I)^{\frac{1}{2}}e^{(s-t)A}x\|
+\int_t^s\|(A+I)^{\frac{1}{2}}e^{(s-\tau)A}b(\tau,x(\tau),u(\tau))\|d\tau
\]
However for every $y\in H$ and $0\leq \tau\leq T$
\[
\|(A+I)^{\frac{1}{2}}e^{\tau A}y\|^2\leq \|(A+I)e^{\tau A}y\| \;
\|y\| \leq \frac{C_1}{\tau^{\delta}}\|y\|^2.
\]
This yields
\[
\|(A+I)^{\frac{1}{2}}e^{\tau A}\|\leq \frac{\sqrt{C_1}}
{\tau^{\frac{\delta}{2}}}
\]
and therefore
\[
\|(A+I)^{\frac{1}{2}}x(s)\|\leq
C_2\left(\frac{1}{(s-t)^{\frac{\delta}{2}}}
+(s-t)^{1-\frac{\delta}{2}}\right)\leq
\frac{C_3}{(s-t)^{\frac{\delta}{2}}}.
\]
We will first show that for every $\epsilon>0$ there exists a
modulus $\sigma_\epsilon$ (also depending on $x$ but independent of
$u(\cdot)$) such that $\|e^{(s_2-s_1)A}x(s_1)-x(s_1)\| \leq
\sigma_\epsilon(s_2-s_1)$ for all $t+\epsilon\le s_1 < s_2\leq T$.
This is now rather obvious since
\[
e^{(s_2-s_1)A}x(s_1)-x(s_1)=\int_0^{s_2-s_1}Ae^{sA}x(s_1)ds
\]
\[
=
\int_0^{s_2-s_1}(A+I)^{\frac{1}{2}}e^{sA}(A+I)^{\frac{1}{2}}x(s_1)ds
-\int_0^{s_2-s_1}e^{sA}x(s_1)ds
\]
and thus
\[
\|e^{(s_2-s_1)A}x(s_1)-x(s_1)\|\leq \|(A+I)^{\frac{1}{2}}x(s_1)\|
\int_0^{s_2-s_1}\frac{\sqrt{C_1}}{s^{\frac{\delta}{2}}}ds
+(s_2-s_1)\|x(s_1)\|
\]
\[
\leq
\frac{C_4}{\epsilon^{\frac{\delta}{2}}}(s_2-s_1)^{1-\frac{\delta}{2}}
+C_5(s_2-s_1).
\]
We also notice that there exists a modulus $\sigma$, depending on
$x$ and independent of $u(\cdot)$, such that
\[
\|x(s)-x\|\le \sigma(s-t).
\]
Let now $t\le s_1 <s_2\leq T$. Denote $\bar s=\max(s_1,t+\epsilon)$.
If $s_2\le t+\epsilon$ then
\[
\|x(s_2)-x(s_1)\|\le 2\sigma(\epsilon).
\]
Otherwise
\begin{multline}
\|x(s_2)-x(s_1)\|\le 2\sigma(\epsilon)+ \|x(s_2)-x(\bar s)\|
\\
\le 2\sigma(\epsilon)+ \|e^{(s_2-\bar s)A}x(s_1)-x(\bar s)\|
+\int_{\bar s}^{s_2}\|e^{(s_2-\tau)A}b(\tau,x(\tau),u(\tau))\|d\tau
\\
\le 2\sigma(\epsilon)+\sigma_\epsilon(s_2-s_1)+C_4(s_2-s_1)
\end{multline}
for some constant $C_4$ independent of $u(\cdot)$. Therefore
(\ref{aaa8}) is satisfied with the modulus
\[
\omega_{t,x}(\tau)=\inf_{0<\epsilon<T-t}
\left\{2\sigma(\epsilon)+\sigma_\epsilon(\tau)+C(\tau)\right\}.
\]

\end{Example}

\end{document}